# Counterexamples in Cake-Cutting

## by Theodore P. Hill

The formal mathematical theory of fair division has a long and rich history dating back at least to Hugo Steinhaus [Ste1 – 4]. Unfortunately, the article "Better Ways to Cut a Cake" [BJK] by Brams, Jones and Klamler in the December 2006 issue of these *Notices* fails to credit the contributions of Steinhaus and other founders of the theory. Worse yet, the article contains many fallacious arguments and conclusions.

The [BJK] article is receiving international media coverage in *Scientific American Science News, Science Daily,* and *the Discovery Channel*, among others, bolstered by our own Society's press release *AMS News Release* [AMS1] and promotional site *AMS in the News* [AMS2]. Other media has given special credibility to the [BJK] paper **because** it was published in the *Notices of the American Mathematical Society*.

Since the *Notices* is so widely read, and represents the Society, it is especially important to set the record straight. Some of the mathematical mistakes in [BJK] can perhaps be salvaged by making simple, but fundamental, changes in both the basic definitions and the hypotheses, and these will be briefly discussed after concrete counterexamples are presented to illuminate some of the errors in the article.

As a basic starting point, consider the classical moving-knife procedure below, attributed in [BJK, p 1319] to Dubins and Spanier, but which in fact is due to Steinhaus's students Stefan Banach and Bronislaw Knaster, as [DS, p2] clearly states: "Their [Banach and Knaster's] solution to the problem is, in essence, as follows."

**Moving-Knife Procedure** [DS, p 2]. "A knife is slowly moved at constant speed parallel to itself over the top of the cake. At each instant the knife is poised so that it could cut a unique slice of the cake. At time goes by the potential slice increases monotonely from nothing until it becomes the entire cake. The first person to indicate satisfaction with the slice then determined by the position of the knife receives that slice and is eliminated from further distribution of the cake. (If two or more participants simultaneously indicate satisfaction with the slice, it is given to any of them.) The process is repeated with the other $n$-1 participants and with what remains of the cake."

The only implicit assumptions here are that the values are nonnegative, they are additive, and, in the direction the knife is moving, they are continuous. The knife need not even be perfectly straight and need not be moved perfectly parallel to itself, and the cake need not be simply connected (as traditional angel-food cakes, with a hole through the center, are not), nor even connected (the cake could have been dropped on the table and landed in pieces).

In the first section of [BJK], the authors state (p 1314) that in problems of fair division of a divisible good, "the well-known 2-person, 1-cut cake-cutting procedure 'I cut, you choose'" is Pareto-optimal, that is, "There is no other allocation that is better for one person and at least as good for the other." Cut-and-choose is not even Pareto optimal among 1-cut procedures, a weaker form of Pareto optimality [BT p 149-150], as the following example shows.

**Counterexample 1**. The "cake" is the unit square, and player 1 values only the top half of the cake and player 2 only the bottom half (and on those portions, the values are uniformly distributed). If player 1 is the cutter, and cuts vertically, his uniquely optimal cut-and-choose solution is to bisect the cake exactly, in which each player receives a portion he values exactly ½. Or if player 1 cuts horizontally, his uniquely optimal risk-adverse cut-and-choose point is the line $y = ¾$, in which case he receives a portion he values at ½ the cake, and player 2 chooses the bottom portion and receives a portion he values at 100% of the cake. But an allocation of the top half of the cake to player 1, and the bottom half to player 2 is at least as good for player 2 in both cases, and is strictly better for player 1, so cut-and-choose is not Pareto optimal in either direction.

In subsequent sections of [BJK], the authors assume that the cake is the unit interval and the value measures are absolutely continuous [BJK p 1315]. But if the cake is indeed the unit interval (which is not the case in classical fair division settings such as [Ste1-4] or [DS], and certainly not the case for the Talmudic scholars who discussed cut-and-choose), the statement [BJK, p 1315] that "We assume that only parallel, vertical cuts, perpendicular to the horizontal x-axis, are made" does not make sense. Even under the hypotheses that the cake is the unit interval (or, equivalently, the unit square with only vertical cuts allowed) and the values are absolutely continuous, the statement in [BJK, footnote 3 p 1318] "an envy-free allocation that uses $n$ -1 parallel cuts is always efficient [i.e., Pareto optimal]", and the corresponding Proposition 7.1 of [BT, p 150], are not true, as the next example shows.

**Counterexample 2**. The cake is the unit interval; player 1 values it uniformly, and player 2 values only the left- and right-most quarters of the interval, and values them equally and uniformly. (In other words, the probability density function (pdf) representing player 1's value is a.s. constant 1 on [0,1], and that of player 2 is a.s. constant 2 on [0, ¼] and on [3/4,1], and zero otherwise.) If player 1 is the cutter, his unique cut point is at $x = ½$, and each player will receive a portion he values at exactly ½. The allocation of the interval [0, ¼] to player 2 and the rest to player 1, however, gives player 1 a portion he values ¾, and player 2 a portion he values ½ again, so cut-and-choose (which is an envy-free allocation for 2 players) is not Pareto optimal.

The two new fair cake-cutting procedures described in [BJK], Surplus Procedure and Equitability Procedure, are not well defined. If a player's value measure does not have a unique median (which absolute continuity does not imply – see Counterexample 2 above), then Surplus Procedure is not well defined. Part (2) of the definition of Equitability Procedure [BJK, p 1318] assumes the existence of "cutpoints that equalize the common value that all players receive for each of the $n!$ possible assignments of pieces to the players from left to right." As the next example shows, such cutpoints may not exist, so Equitability Procedure, too, is not well defined.

**Counterexample 3.** The cake is the unit interval. Player 1 values the cake uniformly, player 2's value is uniform on (0,1/3) (i.e. his pdf is a.s. constant 3 on (0,1/3) and zero elsewhere), and player 3's value is uniform on (2/3,1). Then for the ordering 1-3-2 (from left to right), there do not exist two cutpoints that equalize the values. If the second cutpoint is in (0, 2/3] then player 3 receives 0 but player 1 receives a positive amount. If the second (and hence both cutpoints) are at 0, players 1 and 3 receive 0, and player 2 receives 1. If the second cutpoint is in (2/3, 1), then player 3 receives a positive amount, but player 2 receives 0. If the second cutpoint is at 1, and the first cutpoint is 0, then players 1 and 2 get 0, but player 3 gets 1. If the second cutpoint is at 1 and the first is in (0,1), then player 2 gets 0 but player 1 gets a positive amount. Finally, if both the first and the second cutpoints are at 1, then player 1 gets 1, and both players 2 and 3 get zero.

Even if one imposes extra conditions that guarantee that the system of equations has a solution, step (2) of Eqitability Procedure in [BJK] requires the referee to solve $n$! (possibly highly-nonlinear) systems of $n$-1 integral equations in $n$-1 unknowns. But finding an exact solution, even of *one equation in one unknown,* is not always possible in closed form. For example, if player 1's value is uniform on (0,1) and player 2's value is the standard normal distribution (conditioned to have values in (0,1)), there is no known closed form to the solution of the corresponding integral equation that determines the cutpoint. And without an exact solution, fairness (and equitability and Pareto optimality, etc.) may be lost.

In [BJK, p 1316] the authors define a cake-cutting procedure to be *strategy-vulnerable* if by misrepresenting his value function sent to the referee, a player may "assuredly do better, whatever the value function of the other player", and otherwise the procedure is said to be *strategy-proof*. The article [BJK] contains exactly three new theorems, namely:

**Theorem 1**. *Surplus Procedure is strategy-proof, whereas any procedure that makes e the cut-point is strategy vulnerable.*

**Theorem 2**. *Equitability Procedure is strategy-proof.*

**Theorem 3**. *If a player is truthful under Equitability Procedure, it will receive at least 1/n of the cake regardless of whether or not the other players are truthful; otherwise it may not.*

The following trivial Theorem A shows that the second part of Theorem 1 is false (using Equitability Procedure, which is fair), and that both the first part of Theorem 1 and Theorem 2 are trivial, since both Surplus Procedure and Equitability Procedure are fair.

An allocation procedure is *fair* if each player can guarantee receiving a portion he values at least $1/n$. (An extreme example of an unfair procedure is one that always gives everything to player 1. A procedure which gives the entire cake to each player with probability $1/n$ does give each an *expected* value of $1/n$, but does not *guarantee* any player a portion worth 1/*n* each, so it, too, is not fair in this sense.) It is easy to see that cut-and-choose (for 2 players) and moving-knife are fair procedures.

**Theorem A.** *Every fair procedure is strategy-proof*.

**Proof.** Fix an arbitrary cake-cutting procedure with two players, suppose both players have identical values v, and that both also misrepresent their values identically as a different measure u. Then the procedure allocates a subset S of the cake to one player and allocates its complement ~S to the other. Then at least one player receives a portion he values no more than ½, so that person has not done "assuredly better" than the fair share of ½. Q.E.D.

Even if only one of the players is allowed to misrepresent his strategy, and in case of ties, the pieces are randomly assigned [BJK, (3) p 1315], then every fair procedure is strategy-proof: if player 1 misrepresents his value as u, and u happens to be the true (and declared) value for each of the other players, player 1 will with positive probability receive a portion he values at most 1/*n*, so again he does not do "assuredly better" than fair.

The argument for Theorem 3 is fallacious. The fifth sentence in the proof says that "By moving all players' marks rightward … one can give each player an equal amount greater than 1/n", and the following example shows this is not correct.

**Counterexample 4.** There are three players, the cake is the unit interval [0,1], and all players value it uniformly (i.e., their pdf's are a.s. constant 1 on [0,1]). Then the unique moving-knife marks are at x = 1/3 and x = 2/3, and moving the cuts to the right of those marks will allocate one of the players less than 1/3.

The rest of the argument in [BJK] for Theorem 3 is also incomplete, since the first part only proves a claim about the moving-knife procedure, whereas the desired conclusion concerns Equitability Procedure.

On [BJK, p 1318], the authors claim that their new Equitability Procedure is Pareto optimal (efficient), and on [BJK, p 1320], that their Surplus Procedure is Pareto optimal. Both those claims are false, even when all the value measures are strictly positive everywhere. The underlying reason is that both EP and SP allocate contiguous portions to each player, and as noted in [BT, p 149], "satisfying contiguity may be inconsistent with satisfying efficiency". This is illustrated in the next two examples, which show that EP and SP, respectively, are not in general Pareto optimal.

**Counterexample 5**. The cake is the unit interval [0,1], and there are three players A,B, and C. A's value function is 2.4 on (0, 1/6) and on (1/2, 2/3) and is .3 elsewhere; B's value function is 2.4 on (1/6, 1/3) and (2/3, 5/6) and .3 elsewhere; and C's is 2.4 on (1/3, ½) and (5/6, 1) and .3 elsewhere. Then EP allocates (0, 1/3) to A, (1/3, 2/3) to C, and the rest to B, and each player receives a portion worth exactly .45. But allocating (0, 1/6) and (1/2, 2/3) to A, (1/6, 1/3) and (2/3, 5/6) to B, and the rest to C gives each player a portion he values .

8, which is strictly better for each player than the EP allocation, so EP is not Pareto optimal.

**Counterexample 6**. The cake is the unit interval [0,1], and there are two players A and B. A's value function is 1.6 on (0, 1/4) and on (1/2, 3/4) and is .4 elsewhere; and B's is 1.6 on (1/4, 1/2) and (3/4, 1) and .4 elsewhere. Then SP cuts the cake at 1/2, and each player receives a portion worth exactly .5. But allocating (0, 1/4) and (1/2, 3/4) to A, and the rest to B, gives each player a portion he values .8, which is strictly better for each player than the SP allocation, so SP is not Pareto optimal.

Other arguments in [BJK] are also erroneous or confused. On [BJK, p 1315] the authors "postulate that the players have continuous value functions…and their measures are finitely additive". Does "finitely additive" mean that the Radon-Nikodym theorem may not hold? The authors then state "we assume that the measures of the players are absolutely continuous, so no portion of cake is of positive measure for one player and zero measure for another player". Absolute continuity, of course, even implies *countable* additivity, but the statement that absolute continuity implies that a set of positive measure for one player is also positive for all other players is false, as is easily seen in Counterexample 2 above.

As mentioned above, some of the ideas in the [BJK] article may perhaps be salvaged by making two simple but fundamental modifications: *changing the hypotheses* to require that all value measures (probability density functions on the unit interval) are almost surely strictly positive, and *changing the basic definition* of "strategy-proof" from a strong Pareto-optimality condition to a weaker one.

The difference between strong and weak properties is crucial in many fields – strong/weak topologies, convergence, etc. – and it is in fair division as well, as can be seen here where basic theorems are true using one definition, false using the other. Perhaps replacing the requirement that a strategy of Player A be "assuredly better against all strategies of B", by the weaker requirement that it only be "at least as good against every strategy of B, and strictly better against at least one strategy of B" may lead to interesting correct and non-trivial analogs of Theorems 1 and 2, but even that is not clear. Are the players allowed to bluff only with absolutely continuous measures?

Requiring that the absolutely continuous value measures are also mutually absolutely continuous, or mutually absolutely continuous with respect to Lebesgue measure, may seem like an innocent nonnegativity-versus-positivity technicality, but in fair division such differences often imply important philosophical changes in the problem. For example, suppose that the cake is an inhomogeneous mixture of various ingredients including chocolate. Then absolute continuity simply means that any piece of zero volume is worth zero to every player – a reasonable hypothesis, but not one that is standard in classical or modern fair division theory (e.g., [B, DS, EHK, Hil1-2, J, K, R, RW, Ste1-4, Str]). But further requiring the values to be *mutually* absolutely continuous means that *if one player likes chocolate, all the other players must like chocolate as well*. And requiring that the measures are absolutely continuous with respect to Lebesgue measure (i.e., the corresponding pdf's are almost surely strictly positive), means that *every player must like chocolate, period* (and must like every other part of the cake). It is no great surprise that making such unnatural assumptions may lead to "strong" results.

Perhaps some of the logical errors in [BJK] illuminated above can be corrected by adding additional assumptions, but if the goal is to find "*better* ways to cut a cake", then imposing extra conditions (such as absolute continuity *and* mutual absolute continuity *and* strict positivity of the measures *and* connectivity *and* 1-dimensionality of the cake *and* use of an outside referee), conditions not required by classical procedures like moving-knife or cut-and-choose, can hardly be called an improvement.

The [BJK] article does perhaps contain several new ideas that may serve as a challenge and opportunity for students and mathematicians to develop and make rigorous. As can be seen in the moving-knife procedure above, it is often possible to express practical, yet clean, clear, and beautiful logical conclusions without using highly technical language.

**Acknowledgement.** The author is grateful to Professor Kent E. Morrison and two anonymous colleagues for valuable comments and suggestions.